\documentclass{amsart}
\usepackage{amsmath} 
\usepackage{amsmath}
\usepackage{amsthm}
\usepackage{amssymb}
\usepackage[v2]{xy}

\theoremstyle{plain}
\newtheorem{thm}{Theorem}[section]
\newtheorem{cor}[thm]{Corollary}
\newtheorem{lem}[thm]{Lemma}
\newtheorem{prop}[thm]{Proposition}

\newtheorem*{gq}{General Problem}
\theoremstyle{definition}
\newtheorem{defn}[thm]{Definition}
\newtheorem*{eg}{Example}
\newtheorem*{notation}{Notation}
\theoremstyle{remark}

\newtheorem{obs}[thm]{Observation}

\newcommand{\comment}[1]{}

\newcommand{\G}{\Gamma}

\newcommand{\gammaball}[2]{B_{#1}(#2)}

\newcommand{\gammafinball}[2]{E_{#1}(#2)}

\newcommand{\card}[1]{|#1|}

\newcommand{\ratio}[2]{\mathcal F(#1,#2)}

\newcommand{\genfreeratio}[1]{\mathcal F(#1)}

\newcommand{\coset}[1]{\mathcal \pi^{-1}(#1)}

\newcommand{\idg}{1}

\newcommand{\idc}{I}

\newcommand{\liegpball}[3]{B_{#1}(#2,#3)}

\newcommand{\N}[2]{\mathrm{Nbd}_{#1}(#2)}

\DeclareMathOperator{\Aut}{\mathrm{Aut}}
\DeclareMathOperator{\Isom}{\mathrm{Isom}}
\DeclareMathOperator{\volg}{\mathrm{vol}_G}
\DeclareMathOperator{\volh}{\mathrm{vol}_H}

\newcommand{\Z}{\mathbb Z}
\newcommand{\R}{\mathbb R}

\newcommand{\ra}{\rightarrow}
\newcommand{\rk}{\mathrm{rank}}

\newcommand{\g}{\gamma}

\newcommand{\lcs}[2]{\mathfrak{#1}=\mathfrak {#1}^1 \supset \mathfrak {#1}^2 \supset \cdots \supset \mathfrak {#1}^{#2+1}=0}

\newcommand{\lcsgp}[2]{#1={#1}^1 \supset  {#1}^2 \supset 
\cdots \supset {#1}^{#2+1}=\mathrm{1}_{#1}}

\newcommand{\weight}[1]{\mathcal{W}(#1)}

\newcommand{\polyqb}[1]{Q}
\newcommand{\indpoly}[1]{\psi_{#1}(v)}

\newcommand{\normg}[1]{\|#1\|_{\scriptscriptstyle G}}
\newcommand{\normq}[1]{\|#1\|_{\scriptscriptstyle G/G^k}}
\newcommand{\normgk}[1]{\|#1\|_{\scriptscriptstyle G^k}}
\newcommand{\tatinv}[1]{#1 A(#1^{-1})}
\newcommand{\tatinvmp}{\psi}

\newcommand{\bx}[1]{\mathrm{Box}(#1)}

\newcommand{\num}[1]{\mathcal N(#1)}
\newcommand{\Wa}[1]{\mathcal W_#1}
\newcommand{\W}[2]{\mathcal W_#1^#2}

\begin{document}
\title[The asymptotic density of finite-order elements]{The asymptotic density of finite-order elements in virtually nilpotent groups}

\author{Pallavi Dani}
\begin{abstract}
Let $\G$ be a finitely generated group with a given word metric.    
The \emph{asymptotic density} of elements in $\G$ that have a particular property $P$ is the limit, as $r \ra \infty$, of 
the proportion of elements in the ball of radius $r$ which have the property 
$P$. We obtain a formula to compute the asymptotic density of finite-order elements in any virtually nilpotent group. Further, we show that the spectrum of numbers that occur as such asymptotic densities consists of exactly the rational numbers in $[0,1)$.
\end{abstract}

\maketitle

%

\section{Introduction}
Let $\G$ be a finitely generated infinite group. If $P$ is a
property that elements of $\G$ may have, such as having finite order, having cyclic centraliser or
having a root, it is natural to ask:
What is the density of elements of $\G$ that have the property $P$?

To make this more precise, fix a finite set $S$ of generators for $\G$.
Given two elements $g$ and $h$ in $\G$, set $d(g,h)$ to be the length of the 
shortest word in $S$ representing $g^{-1}h$. This defines the \emph{word metric} on $\G$, which makes $\G$ into a discrete, proper metric space. 
For $r\geq 1$, let $\gammaball{S}{r}$ denote
 the ball of radius $r$ centred at the identity of $\G$ with respect to this
metric. Let $\gammafinball{S}{r}$ denote the set of elements with property $P$ in the ball of
radius $r$.

\begin{gq}
Compute the asymptotics of $\card{\gammafinball{S}{r}}$. 
In particular, find
\[
\ratio{\G}{S}= \lim_{r \rightarrow \infty}
\frac{\card{\gammafinball{S}{r}}}{\card{\gammaball{S}{r}}}
\]
if this limit exists.

\end{gq}
$\ratio{\G}{S}$ is the \emph{asymptotic density} of elements in $\G$ which
have the property $P$. In this paper 
we study the asymptotic density of finite-order elements in the class of
virtually nilpotent groups (i.e. groups containing a nilpotent subgroup of finite index). 
In Theorem \ref{nilp} and Corollary \ref{allvnilpgps}
we obtain a formula to compute $\ratio \G S$ for any virtually nilpotent group $\G$.   

It is worth pointing out that if $\G$ is actually a nilpotent group, the finite-order elements of $\G$ form 
a finite subgroup, so that $\ratio{\G}{S}=0$ for any generating set $S$.
  However, the situation is very 
different when one passes to \emph{virtually} nilpotent groups. 
For example, Theorem \ref{nilp} can be used to show that the densities of finite-order elements in the square and triangle reflection groups in the Euclidean plane are 
$1/4$ and $1/3$, respectively. 
In fact, we prove in Theorem \ref{allrationals} that every rational number in $[0,1)$ occurs as the density of finite-order elements in some virtually nilpotent group. 
This is noteworthy in light of the fact that in many results of this nature in the literature 
the limit is always either $0$ or $1$. A number of such examples are listed in 
\cite{KRSS}. The authors themselves give an example exhibiting 
``intermediate'' density; they show that the union of all proper retracts in 
the free group on two generators has asymptotic density $6/{\pi^2}$.


The phenomenon of positivity of $\ratio \G S$ is not restricted to groups of 
polynomial growth. 
In fact there exist infinite torsion groups with intermediate \cite{Gri} and  even exponential \cite{Adyan,Ly} growth. (For these $\ratio{\G}{S}=1$). 
                
The quantity $\ratio{\G}{S}$ is not a \emph{geometric} property; it may change drastically under quasi-isometry. 
For example, every virtually nilpotent group contains a nilpotent subgroup of finite index (for which $\mathcal F=0$).
However, the 
large-scale geometry of nilpotent Lie groups 
plays an important role in the methods used to study
$\card{\gammafinball{S}{r}}$.
                
  
The idea of studying groups from a
statistical viewpoint was introduced by Gromov, when he indicated that
``almost every'' group is word-hyperbolic. 
Since then the notions of generic group theoretic properties and generic-case behavior have been extensively studied by Arzhantseva, Champetier, 
Kapovich, Myasnikov, Ollivier, Ol'shanskii, Rivin, 
Schupp, Schpilrain, Zuk and others (see \cite{KRSS} and the references contained within).  

\subsection{Virtually nilpotent groups}
Our goal is to compute $\ratio \G S$ for virtually nilpotent groups $\G$. The work of Dekimpe, Igodt, and Lee (See Section \ref{dekigod})
can be used to reduce this question to the following geometric case. 

Let $G$ be a connected, simply connected nilpotent Lie group and 
let $\G$ be a discrete, cocompact subgroup of $\Isom G$, the group of isometries of $G$. 
Auslander \cite{A} generalised Bieberbach's First Theorem to show that $\G$ has a unique maximal normal nilpotent subgroup $\Lambda$, which is torsion-free, and that the quotient $F=\G/\Lambda$ is finite. 


This information determines 
a representation 
$\rho: F \ra \Aut(\mathfrak g)$, where $\mathfrak g$ is the Lie algebra of $G$. 
(See Section \ref{dekigod}). 
If $A$ is an element of $F$, then the automorphism $\rho(A)$ has eigenvalues for its action on $\mathfrak g$. The eigenvalues determined by elements of $F$ in this way depend only on the isomorphism type of $\G$. 
The following theorem gives a formula for $\ratio \G S$ in terms of these eigenvalues. 

\begin{thm}
\label{nilp}
Retaining the above notation,
let $S$ be a finite set of generators for $\G$ 
and let $\gammafinball{S}{r}$ denote the set of finite-order elements in the ball of radius $r$ in the word metric.  
Let
$\lcs {\mathfrak g} k$ be the lower central series of $\mathfrak g$ and let 
$\pi:\G \ra F$ denote the projection map. 

Then there exists $c>0$ such that for any $A\in F$, if $\mathfrak h$ denotes the $1$-eigenspace of $\rho(A)$, then
\begin{equation}\label{cosetbound}
\card{\coset{A} \cap\gammafinball{S}{r}} \leq c r^{d-p}
\end{equation}
where
\[
d=\sum_{i=1}^k i \cdot \rk(\mathfrak g^i/\mathfrak g^{i+1}) \;\; \mathrm{and} \;\;
p=\sum_{i=1}^k i \cdot \rk(\mathfrak h\cap \mathfrak g^i/\mathfrak h \cap \mathfrak g^{i+1})
\]
Further,
\begin{equation}\label{formula}
\ratio{\G}{S}=\frac{m}{\card F}
\end{equation}
where $m$ is the number of elements of $\rho(F)$ that
do not have $1$ as an eigenvalue.
\end{thm}

In particular, $\ratio{\G}{S}$ is independent of the
generating set $S$, so we may write $\genfreeratio{\G}$ instead of $\ratio \G S$.

Every virtually nilpotent group has a subgroup of finite index that acts 
geometrically on a connected, simply connected nilpotent Lie group. (See Section
\ref{dekigod}).  
This is used in Corollary \ref{allvnilpgps} to obtain a formula for $\ratio \G S$ for arbitrary virtually nilpotent groups. 
%

%


The
formula in Theorem \ref{nilp} makes it very easy to compute $\genfreeratio{\G}$ using algebraic data associated with $\G$. \emph{Crystallographic groups}, i.e. groups acting properly discontinuously and cocompactly on Euclidean space,  provide a large class of examples of virtually nilpotent groups. 
There are $17$, $230$, and $4783$ crystallographic groups in dimensions $2$, $3$, and $4$, respectively. 
These are available as libraries  
designed for use with the computer algebra software GAP \cite{gap}. 
The results of the computation of $\genfreeratio{\G}$ for these groups (obtained using GAP) are summarised in the Appendix.

Theorem \ref{nilp} shows that $\genfreeratio{\G}$ is always a rational number. In the following theorem we address the question of which rational numbers in $[0,1]$ can occur.
\begin{thm}
\label{allrationals}
Given any rational number $p/q$ with $0 \leq p/q <1$, there exists a
crystallographic group $\G$ such that 
$\genfreeratio{\G}=p/q$.
\end{thm}
This is proved in Section \ref{examples} by explicitly constructing finite subgroups of $Gl(n,\Z)$ in which exactly $(q-p)/q$ of the elements have eigenvalue $1$. 



The paper is organised as follows.
Sections \ref{background}-\ref{geomnilp} contain definitions and background 
on nilpotent Lie groups. In particular, Section \ref{poly1}
describes certain useful ``polynomial'' coordinate systems for nilpotent Lie groups. Section \ref{poly2} contains some technical lemmas about polynomial coordinates. 

The proof of Theorem \ref{nilp} is contained in Sections \ref{qi}-\ref{finishproof}. 
In Section \ref{qi} we show that a finite-order element of length $r$ in $\G$ fixes a point in a certain ball centered at the identity in $G$. The key 
is to now use the geometry of $G$ to estimate the number of fixed sets of torsion elements that intersect this ball. 
In Section \ref{volumest}, an argument about volumes of balls in $G$ yields 
the upper bound \eqref{cosetbound} for the number of torsion elements in any coset $\coset A$ of 
$\Lambda$. From this bound it follows that if $1$ is an eigenvalue of $\rho(A)$, the torsion in $\coset A$ does not contribute to $\ratio \G S$. 
%
%
%
In Section \ref{finishproof} an inductive argument shows that 
if $1$ is not an eigenvalue of $\rho (A)$, then the coset $\coset A$ consists entirely of torsion elements. 
Theorem \ref{nilp} then follows from the fact that the asymptotic 
density of a coset of $\Lambda$ in $\G$ is $1/\card F$.

Arbitrary virtually nilpotent groups are dealt with in Section \ref{arbit}. Finally, in Section \ref{examples} we construct examples to prove Theorem \ref{allrationals} and also investigate $\genfreeratio \G$ for some non-abelian virtually nilpotent groups. 

This paper consists of a portion of my PhD thesis. I would like to thank my advisor, Benson Farb for his endless support, guidance, and inspiration. 
I would like to thank Angela Barnhill for her suggestions on the manuscript. 


}

\section{Definitions and basic facts}\label{background}
\subsection{Nilpotent Lie groups and Lie algebras}
In this section we recall some background material, which can be found, for example, in \cite{C} or \cite{O}.
Recall that the \emph{lower central series} for a Lie algebra $\mathfrak g$ is defined by
\[
\mathfrak g^1=\mathfrak g, \; \mathfrak g^{i+1}=[\mathfrak g,\mathfrak g^{i}]=
\R\text{-span}\{[X,Y]:X\in\mathfrak g, Y \in \mathfrak g^i\} 
\text{ for } i\geq 1.
\]
Then $\mathfrak g$ is said to be \emph{nilpotent} if $\mathfrak g^{k+1}=\{0\}$ for some $k$. If, in addition, $\mathfrak g^k$ is non-trivial, then $\mathfrak g$ is called a \emph{$k$-step} nilpotent Lie algebra. 

A Lie group $G$ is \emph{nilpotent} if its Lie algebra is nilpotent. 
The lower central series for $G$ is given by $G^1=G$ and $G^{i+1}=[G,G^i]$. 
A connected Lie group is nilpotent if and only if its lower central series is finite. 
If $G^{k+1}=\{1\}$, with 
$G^{k}$ non-trivial, then $G$ is called a \emph{$k$-step} nilpotent Lie group.

A \emph{Lie subgroup} of $G$ is a subgroup which is a submanifold of the underlying manifold of $G$.
If $G$ is connected, the subgroups $G^i$ are Lie subgroups and the Lie algebra of $G^i$ is $\mathfrak g^i$. 
Thus $G$ is $k$-step nilpotent if and only if $\mathfrak g$ is.
For each $i$, the subgroup $G^{i+1}$ is normal in $G^i$ and the quotients $G^{i}/G^{i+1}$ are abelian.

If $G$ is a
connected, simply connected
nilpotent Lie group,
the exponential map, $\exp:\mathfrak g\ra G$,
is an analytic diffeomorphism. Denote its inverse by $\log$.    
Define a map $*:\mathfrak g \times \mathfrak g \ra \mathfrak g$ by
\begin{equation}
X*Y=\log(\exp X \exp Y).
\end{equation}

The \emph{Baker-Campbell-Hausdorff} formula (see \cite{C}, for example) expresses $X*Y$ as a universal power series which involves commutators in $X$ and $Y$. 
While the general term cannot be expressed in closed form,
the low order terms in the formula are well known:
\begin{multline}\label{bakcamphaus}
X*Y=X+Y+\frac 12 [X,Y]+ \frac 1 {12} [X,[X,Y]]-\frac 1 {12}[Y,[X,Y]]\\
-\frac 1 {48} [Y,[X,[X,Y]]] 
- \frac 1{48}[X,Y,[X,Y]]]
+\text{ (commutators in $\geq 5$ terms)}.
\end{multline}
If $G$ is $k$-step nilpotent, then commutators in more than $k$ terms are trivial, which makes this a finite sum.

\subsubsection{Automorphisms and isometries}

\section{Automorphisms and isometries}
}\label{autos-isoms}
An automorphism $A$ of $G$ leaves invariant the groups $G^{i}$. 
Further, $A$ satisfies the relation $A \circ \exp = \exp \circ \,dA$.
The fixed set of $A$ is the image in $G$ of the $1$-eigenspace of $dA$ under the exponential map. It is a Lie subgroup of $G$. 

Let $G$ be endowed 
with a left-invariant Riemannian metric. 
Its group of isometries is
given by $\Isom G= G \rtimes C$, 
where $G$ acts by left multiplication and
$C$ is the group of 
automorphisms of $G$ which preserve the inner product at the identity. 
Thus the action of an element $(g,A)\in \Isom G$ on $t\in G$ is given by 
$(g,A)(t)=gA(t)$.
Any isometry fixing the identity is also an automorphism.

If $G$ is abelian,
then $G=\R^n$ with the standard inner product, where $n$ is the dimension of $G$.  In this case $\Isom G=\R^n \rtimes O(n)$. 

The identity elements of $G$ and $\Aut(G)$ will be denoted by $\idg$ and $\idc$, respectively.  
We will freely make use of the identifications $(g,\idc)=g$ and $(\idg,A)=A$.

Any finite-order isometry of $G$ 
has a fixed point. This follows from a more general result 
of Auslander \cite{ausfixedpt}.
If
 $(g,A)$ is a finite-order isometry with fixed point $p$ (so that $gA(p)=p$), then
\[
(p^{-1},\idc)(g, A)(p,\idc) =(p^{-1}gA(p),A)=(p^{-1}p,A)=(\idg,A).
\]
Thus $(g,A)$ is conjugate to $(\idg,A)$ in $\Isom G$ and hence $\mathrm{Fix}((g,A))=p\mathrm{Fix}(A)$.

\begin{lem}\label{quotient-fixedset}
Let $A=(\idg,A)$ be a finite-order isometry of $G$ fixing the identity. Let $K$ be a normal, $A$-invariant Lie subgroup of $G$ with projection map 
$\pi:G\ra G/K$. If $\bar A$ is the automorphism of $G/K$ induced by $A$, then
$\mathrm{Fix}(\bar A)=\pi(\mathrm{Fix}(A))$. 
\end{lem}

\begin{proof}
Clearly $\pi(\mathrm{Fix}(A)) \subseteq \mathrm{Fix}(\bar A)$. Now if $gK$ is  
fixed by $\bar A$, then $A$ leaves $gK$ invariant. Thus 
$(g^{-1},\idc)(\idg,A)(g,\idc)$ is a finite-order isometry leaving $K$ invariant, which means it has a fixed point in $K$, say $b$. We now have 
\[
(g^{-1}A(g),A)(b)=b \implies g^{-1}A(g)A(b)=b \implies A(gb)=gb.
\]
So $gb$ is fixed by $A$ and
$gK=\pi(gb)$. Thus $\mathrm{Fix}(\bar A)\subseteq \pi(\mathrm{Fix}(A))$.  
\end{proof}


\subsection{Quasi-isometries}
A map $\phi: X\ra X' $ between two metric spaces is a \emph{quasi-isometry} if there exist constants $\lambda \geq 1$, $C\geq 0$, and $D\geq 0$ such that 
\begin{equation*}
\frac 1 \lambda d(x,y) - C \leq d'(\phi(x),\phi(y)) \leq \lambda d(x,y)+C
\end{equation*}
for all $x,y\in X$ and every point of $X'$ is in a $D$-neighbourhood of $\phi(X)$.

The following classical result can be found, for example, in \cite{H}. 
\begin{thm}[Milnor, Efremovich, Schwarzc]\label{fundobs}
If $\G$ is a group acting properly discontinuously and cocompactly by isometries on a proper geodesic metric space $X$, then $\G$ is quasi-isometric to $X$. More precisely, for any $x_0\in X$, the mapping $\G \ra X$ given by $\gamma \mapsto \gamma(x_0)$ is a quasi-isometry. 
\end{thm}

\section{Virtually nilpotent groups}\label{dekigod}
A finitely generated group is said to be
\emph{virtually nilpotent} if it has a nilpotent subgroup of finite
index. 
\emph{Almost crystallographic} groups, i.e.
groups acting properly discontinuously and cocompactly by isometries on a connected, simply connected nilpotent Lie group, are examples of virtually nilpotent groups. This follows from the following theorem of Auslander.

\begin{thm}\label{auslander}\cite{A}
If $\G$ is a discrete, cocompact subgroup of $\Isom G=G\rtimes C$, where $G$ is a connected, simply connected nilpotent Lie group, 
then $\Lambda =\G \cap G$ is cocompact in $G$ and $F=\G/\Lambda$ is a finite group. Further, $\Lambda$ is the unique maximal normal nilpotent
subgroup of $\G$ and it is torsion-free.  
\end{thm}

Using the work of Lee, Raymond, and Kamishima,
Dekimpe and Igodt gave an algebraic condition for a virtually nilpotent group to 
be almost crystallographic.
It is proved in \cite{DI} that 
every virtually nilpotent group 
has a unique maximal finite normal subgroup.  Further, they prove the following:

\begin{thm}\label{leedekig}
If $\G'$ is a virtually nilpotent group with maximal finite normal subgroup $Q$, then $\G=\G'/Q$ is almost crystallographic. 
\end{thm}

This is a generalisation of Malcev's \cite{malcev} result 
that any finitely generated, torsion-free nilpotent group can be embedded as a discrete subgroup of a
nilpotent Lie group, which is unique up to isomorphism. 
Theorem \ref{leedekig} allows us to focus on almost crystallographic groups.

\subsection{Eigenvalues}
Let $\G$ be an almost crystallographic group acting on $G$, i.e. 
there is an injection $\psi:\G \rightarrow G\ltimes \Aut(G)$.
By Theorem \ref{auslander}, $\G$ has a unique maximal normal nilpotent subgroup 
$\Lambda$ with 
$\psi(\Lambda)=G\cap \psi(\G)$, such that $F=\G/\Lambda$ is finite. 

Let $\pi:\G \ra F$ be the projection map. 
There is a unique homomorphism $\xi:F\ra \Aut(G)$ which makes the following diagram commute.

\begin{center}\leavevmode
\xymatrix{
   1 \ar[r]
   &  G \ar[r]
   &  G\ltimes \Aut (G) \ar[r]
   &  \Aut(G) \ar[r]
   &  1  \\
   1  \ar[r]
   & \Lambda \ar[r] \ar[u]^{\psi}
   & \G \ar[r]^\pi \ar[u]^{\psi}
   & F \ar[r] \ar[u] ^{\xi}  
   & 1}
\end{center}      
A diagram-chase
shows that $\xi$ is injective. 
In other words, $F$ can be realised as a group of automorphisms of $G$. 
We obtain an injective homomorphism $\rho:F\ra \Aut(\mathfrak g)$ by composing $\xi$ with the map that assigns to each automorphism in $\Aut(G)$, its 
derivative. 
If $A\in F$, the \emph{eigenvalues} of $A$ are the eigenvalues of the automorphism $\rho(A)$ for its action on $\mathfrak g$.

The fact that these eigenvalues are well-defined follows from a theorem of 
Lee and Raymond \cite{leeray} which says that 
any two isomorphic 
almost crystallographic groups acting on $G$ are conjugate by an element of $G\ltimes \Aut (G)$. 
Indeed if $\psi ':\G\ra G\ltimes \Aut G$ is another injection, giving rise to the homomorphism $\xi ':F \ra \Aut (G)$,
and the element $(g,B)\in G\ltimes \Aut(G)$ conjugates $\psi(\G)$ to 
$\psi '(\G)$, then we  
also have $B\xi(F)B^{-1}=\xi '(F)$, which implies that the eigenvalues assigned to elements of $F$ via $\xi$ are the same as those assigned via $\xi '$.

\section{Polynomial coordinates on $G$}\label{poly1}
For the rest of the paper, $G$ will denote a connected, simply connected nilpotent Lie group.
$G$ can be 
naturally identified with $\R^n$, where $n$ is the dimension of $G$, 
so that the
group structure 
is ``polynomial'' relative to the linear coordinates on $\R^n$.
A map $f:V\ra W$ between two vector spaces is \emph{polynomial} if it is described by polynomials in the coordinates for some (and hence any) pair of bases. 
A \emph{polynomial coordinate map} for $G$ is a diffeomorphism 
$\phi:\R^n\ra G$, 
such that $\log \circ \, \phi$ and $\phi^{-1} \circ \exp$ are polynomial maps. 
We start by defining a useful polynomial coordinate map on $G$. 

Let $\mathfrak g$ be a nilpotent Lie algebra and let $\lcs{g}{k}$ be its lower central series. We define a basis which respects this filtration of $\mathfrak g$.

\begin{defn}
({\bf Triangular basis})
\label{tribasis}
Let $\{X_1, \dots, X_n\}$ be an ordered basis for $\mathfrak g$ with $[X_i,X_j]=\sum_{l=1}^n\alpha_{ijl}X_l$. The basis is \emph{triangular} if $\alpha_{ijl}=0$ when $l\leq \max\{i,j\}$. 
\end{defn}

\begin{eg}
For the three-dimensional Heisenberg Lie algebra (generated by $X, Y$ and $Z$, where $[X,Y]=Z$ and all other brackets are trivial), the ordered sets $\{X, Y, Z\}$ and $\{X+Z, Y, Z\}$ are triangular bases, while the sets $\{ Y, Z, X\}$ and 
$\{X, Y, X+Z\}$ are not.
\end{eg}

A triangular basis can be constructed by starting with an ordered basis for $\mathfrak g^k$ and then 
successively pulling back ordered bases for the factors $\mathfrak g^i/\mathfrak g^{i+1}$, for $i<k$. If $\mathfrak g$ has an inner product, then the triangular basis can be chosen to be orthonormal.  

\begin{defn}({\bf Coordinate map on $G$})\label{polydef}
Let $\{ X_1, \dots, X_n\}$ be a triangular basis for $\mathfrak g$. Define a map $\phi:\R^n\ra G$ by
\[
\phi(s_1, \dots, s_n)=(\exp s_nX_n) \cdots (\exp s_1X_1)= \exp(s_nX_n*\cdots*s_1X_1).
\]
\end{defn}
See \cite[Proposition 1.2.7]{C} for a proof of the fact 
that 
$\phi$ defines a polynomial coordinate map on $G$. 

Each vector $V$ in $\mathfrak g$ is assigned a weight $\mathcal W$, which specifies the smallest group in the lower central series which contains $V$:
\[
\weight V=\max \{i\;|\;V\in \mathfrak g^i\}.
\]

\begin{eg}
In the Heisenberg Lie algebra, with triangular basis $\{X, Y, Z\}$, we have $\weight X=\weight Y=1$ and $\weight Z=2$. 
\end{eg}

In a triangular basis for $\mathfrak g$, there are exactly $\rk {(\,\mathfrak g^i/\mathfrak g^{i+1})}$ vectors which have weight $i$. With this in mind we fix the following notation.

\begin{notation} 
Let $\mathfrak g$ be $k$-step nilpotent and let $\rho_i=\rk{(\,\mathfrak g^i/\mathfrak g^{i+1})}$. A triangular basis for $\mathfrak g$ will be written as $\{X_{ij}\}=\{X_{ij}\;|\;1\leq i\leq k ; 1 \leq j \leq \rho_i\}$, 
where $\weight{X_{ij}}=i$. We will assume that $\{X_{ij}\}$ has the ``dictionary order''. 
Sometimes we will write the basis as $\{X_1, \dots, X_k\}$, where $X_i$ will mean $X_{i1},\dots, X_{i\rho_i}$.

We will identify $G$ with its preimage under the polynomial coordinate map $\phi $ from Definition \ref{polydef}.  Thus the element $s=\exp(s_{k\rho_k}X_{k\rho_k}*\cdots  *s_{11}X_{11})$ of $G$ will be written either as $(s_{ij})$, where it is assumed that $1\leq i \leq k$ and $1 \leq j \leq \rho_i$, or as $(s_1, \dots, s_k)$, where $s_i=s_{i1},\dots,s_{i\rho_i}$ for all $i$.

Finally, we will use $s_i\cdot X_i$ to denote $s_{i\rho_i}X_{i\rho_i}* \cdots* s_{i1}X_{i1}$.
\end{notation}


\section{Geometry of nilpotent Lie groups}\label{geomnilp}

\section{Volumes and distances}\label{metric}

\label{metric}
}

Let $G$ be endowed with a left-invariant Riemannian metric. 
The Ball-Box Theorem of Gromov and Karidi (Theorem \ref{karidimain}) says that in certain polynomial coordinates, the ball of radius $r$ about the identity in $G$
is bounded by certain boxes with sides parallel to the coordinate axes.

Let $\{X_{ij} \;|\; 1\leq i\leq k, 1 \leq j \leq \rho_i\}$ be an orthonormal triangular basis for the Lie algebra $\mathfrak g$ of $G$, where
$\rho_i=\rk{(\mathfrak g^i/\mathfrak g^{i+1})}$. 
Identify $G$ with its preimage under the corresponding 
polynomial coordinate map, and let $\liegpball G \idg r$ denote the ball of radius $r$ about the identity in $G$. 

\begin{defn}
In the above coordinates, for any $l>0$, define $\bx l \subset G$ by
\[
\bx l = \{(s_{ij})\;|\; |s_{ij}|\leq(l)^i \text{ for } 1\leq i \leq k; 1\leq j \leq \rho_i\}.
\]
\end{defn}
This is a box in $G$ with sides parallel to the coordinate axes. For each $i$, it has $\rho_i$ sides of length $2l^i$. Note that the Lebesgue measure of this 
box is $2^n l^d$, where $ n=\sum_{1 \leq i \leq k} \rho_i$ is the dimension of $G$, and $d=\sum_{1\leq i \leq k} i\rho_i$.

\begin{thm}[\bf Ball-Box Comparison Theorem \cite{Gr,K}]\label{karidimain}
There exists $a>1$, which depends only on $G$, 
such that for every $r>1$,  
\begin{equation*} \label{ballbox}
\bx{r/a}
\subset \liegpball G \idg r\\ \subset
\bx{ra}.
\end{equation*}
\end{thm}

The Ball-Box Theorem can be
used to estimate the volume of $\liegpball G \idg r$ and the distances of 
elements of $G$ from the identity. First we make the following definition.

\begin{defn}\label{comparable}
Two functions $f_1$ and $f_2$, from a set $S$ to $\R$ are said to be \emph{comparable}, denoted by $f_1(x)\sim f_2(x)$, 
if there exists $M>1$ such that for all $x \in S$,
\[
\frac 1 M f_2(x) < f_1(x) < M f_2(x).
\]
\end{defn}

There is a unique left-invariant volume form on $G$, up to a scalar multiple. 
Also, the left-invariant measure on $G$ pulls back to
Lebesgue measure on $\R^n$ under the polynomial coordinate map. (See \cite{C}). This yields the following corollary.  

\begin{cor}[\bf Polynomial growth \cite{Gr,K}]\label{polygwth}
Retaining the above notation, if $\volg$ denotes 
the left-invariant volume on $G$, we have
\[
\volg[\liegpball G \idg r] \sim  r^d.
\]
\end{cor}

Let $\normg s$ 
denote the distance of $s\in G$ from  the identity, in the left-invariant metric
on $G$. The following corollary is proved in \cite{ashley}.

\begin{cor}[\bf Distances in nilpotent groups \cite{ashley}]\label{ashley}
Let $s\in G$, with $s=(s_{ij})$ in polynomial coordinates. Then
\[
\normg s \sim \max_{i,j}\{|s_{ij}|^{1/i}\}.
\] 
\end{cor}

If $G^l$ is a group in the lower central series of $G$,
the metric on $G$ induces a left-invariant metric on $G/G^l$. (The inner product at the identity is obtained by identifying $\mathfrak g/\mathfrak g^l$ with 
$\mathfrak g^\perp$).   
The corresponding distances are related as follows: 

\begin{cor}\label{quotientlength}{\bf(Distances in quotients)}
Let $\pi_l:G \ra G/G^l$ be the projection map. 
Then there exists a constant $\delta=\delta(G,l)$, such that for any $s\in G$, 
\[
\|\pi_l(s)\|_{G/G^l} \leq \delta \normg s.
\]
\end{cor}

\begin{proof}
If $\{X_{ij} \;| \;1\leq i\leq k ; 1\leq j \leq \rho_i \}$ 
 is an orthonormal triangular basis
 for $\mathfrak g$ then 
  $\{d\pi_l(X_{ij})\;|\;1\leq i \leq {l-1}; 1\leq j\leq \rho_i\}$
is an orthonormal triangular basis for $\mathfrak g/\mathfrak g^l$,
in the induced left-invariant metric on $G/G^l$. 
 In the corresponding polynomial coordinates, $\pi_l$ is given by $(s_1, \dots,s_k)\mapsto (s_1,\dots ,s_{l-1})$, and the result follows from Corollary \ref{ashley}.
\end{proof}


\section{More on polynomial coordinates}\label{poly2}

}
In this section we 
show that various functions associated with $G$, in particular, the group operations and automorphisms, 
are polynomial maps which 
preserve certain suitably defined weights. (See Proposition \ref{wtpreserving}).
In Lemma \ref{polybound} we obtain a bound on the amount that such weight preserving polynomial maps can stretch distances. 
These results 
are used in the proof of Lemma \ref{fixedptlength}.

We start with an example:
\begin{eg}
Consider the Heisenberg group with polynomial coordinates associated to the triangular basis $\{X, Y, Z\}$. Group multiplication and inversion expressed in these coordinates are given by:
\begin{gather}\label{heismult}
({x}, {y}, {z})
({x_1}, {y_1}, {z_1})
 =({x}+{x_1}, {y}+
 {y_1}, {z}+{z_1}+
 {x}{y_1})\\
\label{heisinv}
({x}, {y}, {z})^{-1}
=(- {x}, - {y}, - {z} +  {x} {y})
\end{gather}
Recall that $\weight X = \weight Y=1$ and $\weight Z=2$. If we assign the weight $1$ to the variables $ {x},  {y}, {x_1}$, and $ {y_1}$ and the weight $2$ to $ {z}$ and $ {z_1}$, then on the right hand side of both 
(\ref{heismult}) and (\ref{heisinv}), the $X$- and $Y$-coordinates are sums of terms of weight $1$, and the $Z$-coordinates are sums of terms such that the total weight of each term is $2$.
\end{eg}

Motivated by this example, we make the following definition.   
Let $y= \{y_{i}\}$ be a set of variables
and let $\mathcal W$ be a function assigning a weight to each $y_i$.
Then polynomials in $\{y_i\}$ can be assigned weights as follows:
\begin{gather*}
\weight {\alpha\, y_{i_1}\cdots y_{i_s}}=
\weight{y_{i_1}}+\cdots+\weight{y_{i_s}}, \text{ where } \alpha \text{ is any constant.}\\
\weight {P(y)} = \max \{\weight{\alpha \,{y_{i_1}\cdots y_{i_s}}}\;|\;\alpha \,y_{i_1}\cdots y_{i_s} \text{ is a term of }P(y)\}.
\end{gather*}

Observe that 
$\weight {P+Q} \leq \max\{\weight P, \weight Q\}$ and 
$\weight{PQ}\leq \weight P + \weight Q$.

\begin{defn}{\bf (Weight-preserving map)}
Let $V$ and $V'$ be vector spaces with bases $\mathcal B=\{X_1,\dots, X_s\}$ and 
$\mathcal B'=\{X'_1, \dots, X'_{s'}\}$ respectively. 
 A polynomial map $f: V \ra V'$ can be written, with respect to these bases, as $f(v)=(P_1(v),\dots , P_{s'}(v))$,
 where $v=(v_1, \dots, v_s) =\sum_{i=1}^s v_iX_i\in V$, and the $P_i$'s are polynomials.
Let $\mathcal W$ (resp. $\mathcal W'$) be a function assigning weights to the $X_i$'s (resp. $X'_i$'s) and define $\mathcal W(v_i)=\mathcal W(X_i)$.  As described above, this induces a weight function $\mathcal W$ on the polynomials $P_i$.
Then $f$ is  \emph{weight-preserving} if $\mathcal W(P_l)\leq \mathcal W'(X'_l)$ 
for all $l$.
\end{defn}
 
\begin{obs}\label{sums-comps}
Finite sums and composites of weight-preserving polynomial maps are weight-preserving polynomial maps. 
\end{obs}

\begin{prop}\label{wtpreserving}
Let $G$ be endowed with a polynomial coordinate system corresponding to the 
triangular basis $\{X_{ij}\;|\;1\leq i \leq k, 1 \leq j \leq \rho_i\}$ of its Lie algebra $\mathfrak g$, where 
$\weight {X_{ij}}=i$. Then the bracket, 
$*$, $\exp$, multiplication and inversion in $G$, and all automorphisms of $G$, when expressed in these coordinates, are 
weight-preserving polynomial maps. 
\end{prop}

\begin{proof}It is easy to see that the bracket is a polynomial map. 
To prove that it is weight-preserving, it is enough to show that 
for given $i$, $j$, $l$, and $m$, if $\alpha$ and $\beta$ are polynomials with 
$\weight \alpha\leq i$ and $\weight \beta \leq l$, then 
$[\alpha X_{ij},\beta X_{lm}]$ is weight-preserving. Observe that 
$[X_{ij},X_{lm}]\in \mathfrak g^{i+l}$, so that  
\[[\alpha X_{ij},\beta X_{lm}]=\sum_{s\geq i+l} \alpha \beta a_{st}X_{st}.
\]
where the $a_{st}$ are structure constants which depend on $X_{ij}$ and $X_{lm}$.   
This is weight-preserving,
since 
$\weight {\alpha \beta a_{st}}\leq i+l \leq s$ for all $s$ and $t$.

Now, the Baker-Campbell-Hausdorff Formula \eqref{bakcamphaus}
expresses $*$ as a finite sum involving brackets. So $*$ is a weight-preserving polynomial map as well, by Observation \ref{sums-comps}.

To prove that $\exp$ is a weight-preserving polynomial map, we produce 
polynomials $Q_{ij}$, with $\weight{Q_{ij}}\leq i$, such that when expressed in
coordinates, 
\begin{equation}\label{expmap}
\exp(v)= (Q_{11}(v), \cdots, Q_{k\rho_k}(v))
\end{equation}
for all $v\in 
\mathfrak g$. 
Recall that $Q_l(v) \cdot X_l$ denotes $Q_{l\rho_l}(v)X_{l\rho_l} * \cdots * Q_{l1}(v)X_{l1}$. 
Define \[
\indpoly l = Q_l(v) \cdot X_l*\cdots * Q_1(v)\cdot X_1
\] 
for all $l$. 
We prove that 
$v=\indpoly k$ (where $\mathfrak g^{k+1}$ is trivial), i.e.
$\exp v= \exp \indpoly k$.
By Definition \ref{polydef}, this is equivalent to equation \eqref{expmap}.
The $Q_{l}$'s are chosen inductively so that 
$\weight{Q_{lj}}\leq l$ and 
$\indpoly l -v \in \mathfrak g^{l+1}$, for all $v$.

Let $v=\sum v_{ij}X_{ij}$ be an element of $\mathfrak g$. 
Set $Q_{1j}(v)=v_{1j}$, for $1 \leq j \leq \rho_1$.
Clearly, $\weight{Q_{1j}}=1$ 
and $\indpoly 1 -v =
\sum_{i>1}v_{ij}X_{ij}\in \mathfrak g^2$. 

Now assume the $Q_{ij}$'s for $i<l$ have been chosen, with
$\indpoly {l-1} -v \in \mathfrak g^l$, say
\begin{equation}\label{indpoly1}
\indpoly {l-1} -v = q_{l1}(v)X_{l1}+ \cdots + q_{l \rho_l}(v)X_{l\rho_l}
+\text{ an element of } \mathfrak g^{l+1}.
\end{equation}
Equivalently, $\indpoly {l-1} -v = (0, \dots, 0, q_{l1}, \dots, q_{l \rho_l}, \dots)$ (this follows from the Baker-Campbell-Hausdorff formula).
%
Observe that $\indpoly {l-1}-v$ 
is a weight-preserving polynomial map, as it is defined in terms of $*$.
Thus $\weight{q_{lj}}\leq l$. 
Choose $Q_{lj}= -q_{lj}$ for $1 \leq j \leq \rho_l$. Then 
$\weight{Q_{lj}} \leq l$.  Moreover, 
using the
Baker-Campbell-Hausdorff formula again, we have
\begin{equation}\label{indpoly2}
\begin{aligned}
\indpoly l 
&= Q_{l\rho_l}(v)\,X_{l\rho_l}*\cdots \cdot *Q_{l1}(v)\,X_{l1}*\indpoly {l-1}\\
&= Q_{l1}(v)\,X_{l1}+\cdots \cdot +Q_{l \rho_l}(v)\,X_{l\rho_l}+\indpoly {l-1}+
\text{ an element of } \mathfrak g^{l+1}\\
&= -q_{l1}(v)\,X_{l1}+\cdots \cdot -q_{l \rho_l}(v)\,X_{l\rho_l}+\indpoly {l-1}+
\text{ an element of } \mathfrak g^{l+1}.
\end{aligned}          
\end{equation}
Equations \eqref{indpoly1} and \eqref{indpoly2} imply that 
$\indpoly l -v \in \mathfrak g^{l+1}$, completing the induction. Since 
$\mathfrak g^{k+1}$ is trivial, we have $v=\indpoly k$ as required.

%

Now let $s=(s_{ij})$
and $t=(t_{ij})\in G$. Using Definition \ref{polydef},
multiplication and inversion can be written in terms of $*$ and $\exp$ as follows:
\begin{align*}
st&=
\exp (s_k \cdot X_k* \cdots* s_1 \cdot X_1*
t_k \cdot X_k* \cdots* t_1 \cdot X_1)\\
s^{-1}&=\exp[-(s_k\cdot X_k*\cdots*s_1\cdot X_1)]
\end{align*}
%
If $A$ is an automorphism,
$dA$ preserves the bracket, and hence $*$. 
Thus  
\begin{align*}
A(s)
&=A(\exp(s_{k\rho_k}X_{k\rho_k}*\cdot \cdot *s_{k1}X_{k1}*\cdots*s_{1\rho_1}X_{1\rho_1}*\cdot \cdot *s_{11}X_{11}))\\
&=\exp dA(s_{k\rho_k}X_{k\rho_k}*\cdot \cdot *s_{k1}X_{k1}*\cdots*s_{1\rho_1}X_{1\rho_1}*\cdot \cdot *s_{11}X_{11})\\
&=\exp (s_{k\rho_k}dAX_{k\rho_k}*\cdot \cdot *s_{k1}dAX_{k1}*\cdots*s_{1\rho_1}dAX_{1\rho_1}*\cdot \cdot *s_{11}dAX_{11}).
\end{align*}
Now for each $i$ and $j$, we have $dA(X_{ij}) \in \mathfrak g^i$, so that
\[
s_{ij}dAX_{ij}= \sum_{\substack {l\geq i\\ 1\leq m\leq\rho_l}} s_{ij}\alpha_{lm}X_{lm}.
\]
where the $\alpha_{lm}$'s are constants depending on $A$. Thus $s\mapsto s_{ij}dAX_{ij}$ is a weight-preserving polynomial map.

It follows from Observation \ref{sums-comps} that multiplication inversion and
all automorphisms are weight-preserving polynomial maps. 
\end{proof}

We now obtain a bound on the amount that a weight-preserving polynomial map can stretch distances.

\begin{lem}\label{polybound} 
If $P:G \ra G$ is a weight-preserving polynomial map, then
there exists a constant $\lambda=\lambda(P)>0$ such that for all $y \in G$,
\[
\normg {P(y)} \leq \lambda \normg y.
\]
\end{lem}

\begin{proof}
By Corollary \ref{ashley}
we know that there exists $\mu >1$ such that 
\[
\frac 1 \mu \,\normg y \leq \max_{i,j}\{|y_{ij}|^{1/i}\} \leq \mu \,\normg y
\;\;\;\text{ and }\;\;\;
\frac 1 \mu \,\normg P \leq \max_{i,j}\{|P_{ij}(y)|^{1/i}\} \leq \mu \,\normg P.
\]
Let $\alpha y_{i_1j_1}\cdots y_{i_lj_l}$ be a term occurring in $P_{ij}$ for some $i$ and $j$. We omit the second subscript for convenience. The weight-preserving condition, 
$\weight{P_{ij}}\leq i$, implies that $i_1+\cdots+i_l\leq i$. 
Let $s$ be such that 
$
|y_{i_s}|^{1/{i_s}}=\max \{ |y_{i_1}|^{1/{i_1}},\cdots, |y_{i_l}|^{1/i_l}\}$. 
Then 
\begin{align*}
|\alpha y_{i_1}\cdots y_{i_l}|^{1/i} &=|\alpha|^{1/i} \left[\left(|y_{i_1}|^{1/i_1}\right)^{i_1}\cdots \left(|y_{i_l}|^{1/i_l}\right)^{i_l}\right]^{1/i}\\
&\leq |\alpha|^{1/i}
\left[|y_{i_s}|^{i_1/{i_s}}\cdots|y_{i_s}|^{i_l/{i_s}}\right]^{1/i}\\
&= |\alpha|^{1/i} \left[|y_{i_s}|^{(i_1+\cdots+i_l)/{i_s}}\right]^{1/i}
\leq |\alpha|^{1/i} |y_{i_s}|^{1/{i_s}}
\leq |\alpha|^{1/i} \mu \normg y.
\end{align*}
This, combined with the fact that  
$
|P_{ij}(y)|^{1/i} \leq \sum_{\text {terms of }P_{ij}(y)}
|\alpha y_{i_1}\cdots y_{i_l}|^{1/i}$, 
enables us to choose a constant $\nu$ such that 
$|P_{ij}(y)|^{1/i}\leq \nu \mu \,\normg y$
for all $i$ and $j$. Thus $\normg P \leq \mu^2 \nu \normg y$, 
and we can take $\lambda=\mu^2\nu$. 
\end{proof}



\section{Discrete subgroups of $\Isom G$}
}
\section{The relation between $\G$ and $G$}\label{qi}
We now return to the set-up in Theorem \ref{nilp}. 
Let
$\G$ be a discrete, cocompact subgroup of $\Isom G=G \rtimes C$, with maximal 
normal nilpotent subgroup $\Lambda$ and finite quotient $F=\G/\Lambda$, which we 
identify with its image in $\Aut (G)$ under the map $\xi$ defined in Section 
\ref{dekigod}.   

Let $S$ be a 
finite generating set for $\G$ and let $\gammafinball{S}{r}$ be the set of finite-order elements of length less than or equal to $r$ in the word metric. 

Every finite order isometry of $G$ has a fixed point (see Section \ref{autos-isoms}). 
The following lemma will allow us to estimate the 
cardinality of
$\gammafinball{S}{r}$ by counting fixed sets of elements in $\gammafinball{S}{r}$ that intersect a certain ball in $G$.

\begin{lem}\label{qilemma}
Let $G$, $\G$, $S$ and $\gammafinball{S}{r}$ be as above. Then 
there exists a positive constant
$\kappa$ such that 
if 
$(g,A) \in \gammafinball{S}{r}$, then the fixed set of $(g,A)$ in $G$ intersects 
$\liegpball G \idg {\kappa r}$.
\end{lem}

The proof relies on the following proposition, which establishes a 
relation between the displacement of a point under the action of a finite-order
isometry fixing the identity of $G$, and the distance of that point from the fixed set of the isometry.


\begin{prop}\label{fixedptlength}
Let $A$ be a finite-order isometry of $G$ fixing the identity, with fixed subgroup $H$. Then there exists $K=K(A,G)$, such that
for all $r>0$,  
if $t \in G$ satisfies  
 $\normg{\tatinv t}<r$, then
 there exists $h\in H$ such that
  $\normg{th}<Kr$

\end{prop}

\begin{proof}[Proof of Proposition \ref{fixedptlength}]
The proof is by induction on the lower central series of $G$.  Firstly, if $G$ is abelian, we may  
 write the condition on $t$ as 
$\normg{t-A(t)}<r$, where $A\in O(n)$. 
In this case, let $H^\perp$ be the orthogonal complement of $H$ in $G$. 
There exists $h\in H$ such that $t+h \in H^\perp$. Since $h$ is a fixed point of $A$, we have 
\begin{equation}\label{t+h}
\normg{(t+h)-A(t+h)}=\normg{t-A(t)}<r.
\end{equation}

The map $A$ 
leaves $H^\perp$ invariant and has no fixed points on the compact set $\{x\in H^\perp|\,\normg{x}=1\}$. Thus the function $\normg{x-A(x)}$ attains a positive minimum, say $m$, on this set. 
Now,
\[  
\normg{(t+h)-A(t+h)}=\normg{t+h} \left\| \frac{t+h}{\|{t+h}\|\rlap{$_{\scriptscriptstyle G}$}}-A\left(\frac{t+h}{\|{t+h}\|\rlap{$_{\scriptscriptstyle G}$}}\;\,\right)\right\|_G  \geq m \normg{t+h}.
\]
Inequality \eqref{t+h} now implies that $\normg{t+h}\leq \left(\frac 1 m \right) r $.

Now let $G$ be $k$-step nilpotent 
with an orthonormal 
polynomial coordinate system as in Section \ref{ballbox}.
Let $\pi:G\ra G/G^k$ be the canonical projection, i.e. $\pi(x_1, \dots x_k)=(x_1,\dots, x_{k-1})$, and let $\bar A$ be the automorphism of $G/G^k$ induced by $A$. 

Let $t\in G$ with $\normg{\tatinv t}<r$.
For the rest of the proof, we use 
$\prec$ to mean \emph{less than, up to a constant factor that depends only on $G$ and $A$}. We produce $h\in H$ such that $\normg{th} \prec r$. 

Corollary \ref{quotientlength} on distances in quotients, implies that
\[
\normq {\pi (t) \bar A\left(\pi(t)^{-1}\right)} =\normq{\pi(\tatinv t)  } 
 \prec \normg{\tatinv t} < r.
\]
By Lemma \ref{quotient-fixedset}
the fixed set of $\bar A$ is $\pi(H)$. 
So by the induction hypothesis, there exists 
$h_1\in H$ such that 
$\normq{\pi(t)\pi(h_1)}\prec r$.

We may write $th_1=yz$, where $y=(y_1, \dots, y_{k-1},0)$ and $z=(0, \dots,0,z)\in G^k$. 
Note that $\pi(y)=\pi(th_1)=\pi(t)\pi(h_1)$, so that $\normq{\pi(y)}\prec r$.
Further, Corollary \ref{ashley} implies that
\[ \normg y  \sim \normq{\pi(y)} \prec r.
\]
However, $\normg z$, and hence $\normg{th_1}$, may be arbitrarily large. This will be fixed by correcting $th_1$ by an element of 
 $H \cap G^k$.

We first show that $\normg{\tatinv z} \prec r$.
Note that $z$ is in the centre of $G$, which is preserved by $A$, so that 
\begin{equation}\label{yz}
yz A\left([yz]^{-1}\right)= yzA(z^{-1}y^{-1})=\tatinv y \tatinv z.
\end{equation}

Proposition \ref{wtpreserving} and Observation \ref{sums-comps}
imply that $x \mapsto \tatinv x$ is a 
 weight-preserving polynomial map. 
Lemma \ref{polybound} now implies that
$\normg {\tatinv y} \prec \normg y \prec r$.
Further, $yz A\left([yz]^{-1}\right)=
th_1A\left([th_1]^{-1}\right)=
\tatinv t$, so that $\normg {yz A\left([yz]^{-1}\right)} \prec r$. 
Equation \eqref{yz} now implies that  
$\normg{\tatinv z} \prec r$. 

Corollary \ref{ashley} implies that 
$\normgk{x} \sim \normg{x}^k$ for all $x \in G^k$.
Thus
 $\normgk{\tatinv z} \prec r^k$. 
By the first step of the induction, (since $G^k$ is abelian) there exists $h_2 \in H \cap G^k$, such that $\normgk{zh_2}\prec r^k$, which means 
$\normg{zh_2}\prec r$. 
Setting $h=h_1h_2$ completes the inductive step:
\begin{equation*}
\normg{th}=
\normg{th_1h_2} = \normg{yzh_2} \leq \normg{y}+\normg{zh_2}\\ 
\prec r.
\end{equation*}
\end{proof}
 

\begin{proof}[Proof of Lemma \ref{qilemma}]
Let $\ell_S$ denote the distance from the identity in $\G$. 
Since the map $\G \ra G$ given by $\gamma \mapsto \gamma(\idg)$ is a 
quasi-isometry
(Theorem \ref{fundobs}), there exist positive constants $\lambda$ and $C$, such that for all $(g,A)\in \G$,
\begin{equation*}\label{quasiisom}
\frac 1 \lambda \ell_S((g,A)) - C \leq \normg g \leq \lambda \ell_S((g,A)) +C.
\end{equation*}
So if $(g,A)$ is an element of $\gammafinball{S}{r}$, then 
$\normg g \leq \lambda r+ C \leq (\lambda +C)r$ (since $r \geq 1$).

Let $H$ be the fixed subgroup of $A$. As discussed in Section \ref{autos-isoms},  
if $(g,A)$ fixes $t_g$, then its fixed set is $t_gH$,
and $g=t_gA(t_g^{-1})$. Thus $\normg{\tatinv{t_g}}<(\lambda+C)r$.
So by Proposition \ref{fixedptlength} there exists $K=K(A,G)$,
and an element $h\in H$ such that 
\[
\normg{t_gh}<K(\lambda + C)r.
\] 
In other words, the fixed set of $(g,A)$ intersects $\liegpball G \idg {\kappa r}$, where 
$\kappa=K(\lambda+C)$. 
Since $F$ is finite, 
$K$ can be chosen to work simultaneously for all $A\in F$.
\end{proof}


\section{Volume estimates in $G$}\label{volumest}

\section{Volume estimates in $G$}

}
In this section we fix an element $A$ in the finite quotient $F$, and obtain an
upper bound for $\card {\coset A \cap \gammafinball{S}{r}}$. 
Let
$\mathcal H$ be the collection of fixed sets in $G$ of finite-order elements in $\coset A$. Then Lemma \ref{qilemma} says that every 
element of $\mathcal H$ intersects $\liegpball G \idg R$, where $R=\kappa r$.

If $H$ is the fixed subgroup of $A$, then $\mathcal H$ consists of 
cosets of $H$ in $G$. 
We 
will use a volume argument to count the number of elements of $\mathcal H$ intersecting $\liegpball G \idg R$.  
We first obtain disjoint neighbourhoods of the submanifolds in $\mathcal{H}$ and intersect them with 
$\liegpball{G}{\idg}{R}$.
Then we use the fact that the volume of 
$\liegpball{G}{\idg}{R}$
is greater than the sum of the volumes of these disjoint pieces contained in it.

\begin{lem}\label{disjoint}
Let $H$ and $\mathcal H$ be as above. Then
there exists $\epsilon>0$ such that the $\epsilon$-neighbourhoods of
cosets in $\mathcal{H}$ are pairwise disjoint. 
\end{lem}

\begin{proof}
We first show that 
if $t \notin H$, then 
$H$ and $tH$ have disjoint $\delta$-neighbourhoods for some $\delta>0$.  
If there is no such $\delta$, we can find sequences $\{h_i\}$ and $\{tk_i\}$, with $h_i,k_i \in H$, such that 
$d(h_i, tk_i)\rightarrow 0$, which
  means $\{h_i^{-1}tk_i\}$
  is a sequence converging to $\idg$. We now have
\begin{align*}
A(h_i^{-1}tk_i) \rightarrow A(\idg) 
&\implies h_i^{-1}A(t)k_i \rightarrow \idg\\
&\implies k_i^{-1}A(t^{-1})h_i \rightarrow \idg\\
&\implies (h_i^{-1}tk_i)(k_i^{-1}A(t^{-1})h_i) \rightarrow \idg\\
& \implies h_i^{-1}tA(t^{-1})h_i \rightarrow \idg. 
\end{align*}
Set $g=tA(t^{-1})$. Suppose $g\in G^j$.
Then
$\{h_i^{-1}gh_i\}=\{g(g^{-1}h_i^{-1}gh_i)\}$ is a sequence in $gG^{j+1}$
converging to $\idg$.
Since $gG^{j+1}$ is closed set, the limit, $\idg$, is in $gG^{j+1}$. Thus $g \in G^{j+1}$. This inductive argument shows that  
$g=\idg$. 
This means that $t$ is a fixed point of $A$, contradicting the fact that $t \notin H$.

Since $\mathcal H$ is a discrete collection of cosets of $H$, the above implies the existence of $\epsilon>0$ such that for any $tH \in \mathcal H$, the $\epsilon$-neighbourhoods of $H$ and $tH$ are disjoint. It follows that the $\epsilon$-neighbourhoods of any two cosets in $\mathcal H$ are disjoint. 
\end{proof}

We denote the $\epsilon$-neighbourhood of a set $Y$ by $\N \epsilon Y$.
We will need to estimate the volumes of intersections of $\epsilon$-neighbourhoods of elements of $\mathcal H$ with $\liegpball G \idg R$. 
Since $\N \epsilon {tH}=t \N \epsilon H$, we will focus on estimating the volume of $\N \epsilon H \cap \liegpball G \idg R$. 

\comment{
The left-invariant metric on $G$ induces a left-invariant metric on the Lie subgroup $H$.  Note that the distance between two points in $H$ measured in the $G$-metric may be less then their distance in the induced metric on $H$. If $h\in H$ and $\liegpball H h R$ denotes balls of radius $R$
in the $H$-metric, then $ \liegpball H 1 R$ is, in general, a subset of $\liegpball G h R \cap H$. 
}

The following lemma
relates the volume of $\N \epsilon H \cap \liegpball G \idg R$ to the volume 
of 
$H \cap \liegpball G \idg R$ with respect to the left-invariant measure on $H$. 

\begin{lem}
\label{intest}
Let $\volg$ and $\volh$ denote the left-invariant volumes in $G$ and $H$, respectively. Let $\epsilon$ be the constant obtained in Lemma \ref{disjoint}.
There exists a constant $V_\epsilon>0$, which is independent of $R$, such that 
\[
\volg\left[\N{\epsilon}{H} \cap \liegpball{G}{\idg}{R}\right]>V_{\epsilon}  \volh\left[H\cap\liegpball{G}{\idg}{R}\right].
\]
\end{lem}

\begin{proof}
For any $R$, there exists a finite set $\mathcal D_R$ of points in $H \cap \liegpball{G}{\idg}{R-\epsilon}$ which satisfies the following two conditions: 
\begin{enumerate}
\item \label{smallballs}
Balls of radius $\epsilon$ in $G$, centred at points in $\mathcal D_R$ are disjoint. 
\item \label{3epsilon}
Balls of radius $3\epsilon$ in $G$, centred at points in $\mathcal D_R$ cover $H\cap \liegpball{G}{\idg}{R}$.
\end{enumerate}

Note that each $\epsilon$-ball as in (\ref{smallballs}) is contained in $\N{\epsilon}{H}\cap \liegpball{G}{\idg}{R}$ and has volume equal to $V_1^\epsilon=\volg[\liegpball{G}{\idg}{\epsilon}]$.
Thus we have 
\begin{equation}\label{ineq}
\volg[\N{\epsilon}{H}\cap \liegpball{G}{\idg}{R}] > V_1^{\epsilon} \;\card{\mathcal D_R}.
\end{equation}
If $h\in H$, then $\liegpball G h {3\epsilon} \cap H = h \left(\liegpball G \idg {3\epsilon} \cap H\right)$.  
Thus, the volume in $H$ of 
the intersection of $H$ with
 a $3\epsilon$-ball as in (\ref{3epsilon}) 
is a constant,
$V_2^\epsilon=\volh[\liegpball G \idg {3\epsilon} \cap H]$.
Since the collection of balls in (\ref{3epsilon}) cover 
$H \cap \liegpball{G}{\idg}{R}$, we have
\begin{equation}\label{lowerbound}
V_2^\epsilon \;
\card{\mathcal D_R}>\volh[H \cap \liegpball{G}{\idg}{R}].
\end{equation}

Set $V_{\epsilon}=\frac {V_1^\epsilon}{V_2^\epsilon}$. Combining inequalities (\ref{ineq}) and (\ref{lowerbound}) yields the result:
\[
\volg[\N{\epsilon}{H}\cap \liegpball{G}{\idg}{R}] >
V_\epsilon
\volh[H \cap \liegpball{G}{\idg}{R}].
\]
\end{proof}

The next step is to estimate $\volh [H\cap\liegpball{G}{\idg}{R}]$. 
Note that the distance between two points in $H$ measured in the metric on $G$ may be less then their distance in the induced metric on $H$. If $\liegpball H \idg R$ denotes the ball of radius $R$ in the 
induced metric on $H$, then $ \liegpball H \idg R$ is, in general, a subset of $H \cap \liegpball G \idg R$. 

\subsection{Polynomial coordinates compatible with $H$}
We will define a new polynomial coordinate system on $G$, such that the preimage of $H$ under the polynomial coordinate map is a subspace of $\R^n$ parallel to the coordinate axes. We will then be able to use the ball-box technique from Theorem \ref{karidimain} to estimate volumes in $H$ and $G$ simultaneously.

Let $\mathfrak h$ be the Lie algebra of $H$. We
choose a triangular basis for $\mathfrak g$, such that a subset of the basis is a triangular basis for $\mathfrak h$. This can be done as follows. 

Let $\lcs{g}{k}$ be the lower central series of $\mathfrak g$. Let  \[\rho_i=\rk(\mathfrak g^i/\mathfrak g^{i+1}) \text{ and } \eta_i=\rk(\mathfrak h\cap \mathfrak g^i/\mathfrak h \cap\mathfrak g^{i+1}).\]
For each $i$, pick $X_{i1}, \dots,X_{i\rho_i}$ to be a pullback of a basis for $\mathfrak g^i/\mathfrak g^{i+1}$, such that $X_{i1},\dots, X_{i\eta_i}$
projects to a basis for $\mathfrak h\cap \mathfrak g^i/\mathfrak h \cap\mathfrak g^{i+1}$. Give $\{X_{ij}\;|\;1\leq i \leq k; 1\leq j\leq \rho_i\}$ the dictionary order. 

It is easy to see that this gives a triangular basis for $\mathfrak g$.
Since $H$ is a Lie subgroup, $\mathfrak h$ is a subalgebra. In particular, it is closed under the bracket, so that $\{X_{ij}\;|\;1\leq i \leq k; 1\leq j\leq \eta_i\}$
is a triangular basis for $\mathfrak h$. 

Now define a polynomial coordinate map $\phi:\R^n \ra G$ as in Definition
\ref{polydef}. Observe that $\phi^{-1}(H)$ is the set of points $\{(s_{ij})\in \R^n\;|\; s_{ij}=0 \text{ if } \eta_i<j \leq \rho_i \}$, which is a plane spanned by a subset of the coordinate axes for $G$.

We now endow $G$ with a new left-invariant metric that makes the above basis orthonormal. Note that, up to a constant factor, there is only one left-invariant volume form on a Lie group. Since we are only interested in the degree of growth, we may use this new metric to estimate volume. 

Recall that the symbol $\sim$ denotes comparable functions. (See Definition \ref{comparable}).

\begin{lem}
\label{hvolume}
Retaining the above notation, let $\displaystyle p=\sum_{i=1}^k i \; \eta_i$. Then
\[
\volh[H\cap\liegpball{G}{\idg}{R}] \sim R^p.
\] 
\end{lem}

\begin{proof} Theorem \ref{karidimain} tells us that in the coordinate system defined above, $\liegpball G\idg R$ can be bounded by two boxes (one contained in it and one containing it) which have sides parallel to the coordinate axes. In particular,
there exists $a>1$ such that for $R>1$,
\begin{equation*} 
\{ (s_{ij})\;|\; |s_{ij}|\leq(R/a)^i \text{ for all } i\} 
\subset \liegpball G \idg R\\ \subset
\{ (s_{ij})\;| \;|s_{ij}|\leq(aR)^i \text{ for all } i\}. 
\end{equation*}
For each $i$, the outer box has $\rho_i$ sides of length $2(aR)^i$. The intersection of this box with $H$ is a box parallel to the coordinate axes in $H$, with $\eta_i$ sides of length $2(aR)^i$, for each $i$. The Lebesgue measure of this intersection is therefore a constant multiple of $R^p$, where $ p=\sum_{i=1}^k i \, \eta_i$.
A similar statement holds for the inner box. 
Moreover,
$H\cap \liegpball G \idg R$ is contained in the outer box and contains the inner box. 
This proves the lemma, since the Lebesgue measure on $\phi^{-1}(H)$ is comparable to the left-invariant measure on $H$.
\end{proof}

We can now prove inequality \eqref{cosetbound} in Theorem \ref{nilp}.

\begin{lem}\label{estimate}
There exists $c>0$ such that 
$
\card {\coset A \cap \gammafinball{S}{r}} \leq c r^{d-p}
$,
where $d=\sum_{i=1}^k i \; \rho_i$ and $p=\sum_{i=1}^k i \; \eta_i$.
\end{lem}

\begin{proof}
In Lemma \ref{disjoint} we obtained disjoint $\epsilon$-neighbourhoods of the cosets in $\mathcal H$. For every $tH \in \mathcal H$ 
which intersects $\liegpball G \idg R$,
choose an element $p_t$ in the intersection. Observe that 
$\N{\epsilon}{tH}=p_t(\N{\epsilon}{H})$, so that the sets $p_t(\N{\epsilon}{H}\cap \liegpball{G}{\idg}{R})$ corresponding to distinct cosets in $\mathcal H$ are disjoint. 
Since $\normg{p_t}\leq R$, we have
\[
p_t(\N{\epsilon}{H}\cap \liegpball{G}{\idg}{R}) \subseteq \liegpball{G}{\idg}{2R}.
\]
Let $M$ be the number of elements of $\mathcal{H}$ intersecting $\liegpball{G}{\idg}{R}$. Then 
\begin{align*}
\volg[\liegpball{G}{\idg}{2R}]
& > \volg \left[ \bigcup_
{t\mathcal H \cap \liegpball G \idg R \neq \emptyset}
p_t(\N{\epsilon}{H}\cap \liegpball{G}{\idg}{R})\right]\\
& =M \volg[\N{\epsilon}{H}\cap \liegpball{G}{\idg}{R}]\\
& >M V_{\epsilon}\volh[H\cap\liegpball{G}{\idg}{R}], & (\text{Lemma} \;\ref{intest})
\end{align*}
so that 
\[
M < \left(\frac 1 {V_\epsilon}\right)
\frac{\volg[\liegpball{G}{\idg}{2R}]}{\volh[H\cap\liegpball{G}{\idg}{R}]}.
\]
Corollary \ref{polygwth}
and Lemma \ref{intest} now imply the existence of a constant $c'>0$ such that 
$M<c' R^{d-p}$. Now by Lemma \ref{qilemma}, $M$ is an upper bound for 
$\card {\coset A \cap \gammafinball{S}{r}}$, so that 
$\card {\coset A \cap \gammafinball{S}{r}} \leq cr^{d-p}$, where $c=\kappa c'$.
\end{proof}


\section{Finishing the proof}\label{finishproof}

\section{Proofs of Theorems blah and blah}

}
To complete the proof we will need the following theorem of Pansu on the growth of balls in virtually nilpotent groups. 

\begin{thm}\cite{P}\label{pansu}
Let $\G$ be a finitely generated, virtually nilpotent group with finite generating set $S$. Let $d=\sum_{i=1}^\infty i\, \rk(\G^i/\G^{i+1})$.
Then $\lim_{r\ra \infty} \frac{\card{\gammaball{S}{r}}}{r^d}$ exists. 
\end{thm}

In particular, $\card{\gammaball S r} \sim r^d$.
Together with 
Lemma \ref{estimate}, this implies that if $1$ is an eigenvalue of $dA$, then 
\begin{equation}\label{contribution}
\lim_{r \ra \infty} \frac{\card{\coset A \cap \gammafinball{S}{r}}}{\card {\gammaball{S}{r}}} = 0.
\end{equation}
This is because in this case, the fixed set $H$ of $A$ has dimension at least 
$1$, which implies that $p\geq 1$ and hence $d-p<d$.

On the other hand if $1$ is not an eigenvalue of $dA$, we have the following. 
\begin{lem}\label{fullcoset}
Let $A$ be a finite-order isometry fixing the identity, such that 
$1$ is not
an eigenvalue of $dA$. Then $(g,A)$ has finite order for every $g \in G$.  
\end{lem}

\begin{proof}
We prove in Lemma \ref{surj} below, that for every $g \in G$, there exists 
$t\in G$ with $g=\tatinv t$. Now 
$(t,\idc)(\idg,A)(t^{-1},\idc)=(tA(t^{-1}),A)=(g,A)$. In other words, $(g,A)$
is conjugate in $\Isom G$ to $(\idg, A)$, and hence has finite order. 
\end{proof}

\begin{lem}\label{surj}
Let $A$ be a finite-order isometry fixing the identity, such that $1$ is not an
eigenvalue of $dA$. 
The map $\tatinvmp:G \rightarrow G$ defined by $\tatinvmp(t)=tA(t^{-1})$
is surjective. 
\end{lem}

\begin{proof} The proof is by induction on the lower central series. 
If $G$ is abelian, we can
write $\tatinvmp(t)=t-A(t)=(I-A)t$,
where $A=dA$ is linear. Since $1$ is not an eigenvalue of $A$,
there is no non-zero $v$ with $(I-A)v=0$. Thus $I-A$ is invertible and the equation $\tatinvmp(t)=b$ has a solution for every $b$.

Now let $\lcsgp{G}{k}$ be the lower central series for $G$. 
Then $\tatinvmp$ leaves $G^i$ invariant for all $i$, since $A$ does. 
Assume $\tatinvmp|_{G^i}:G^i\ra G^i$ is surjective.

The automorphism $A$ induces an 
automorphism $A_i$ on $G^{i-1}/G^i$.
It follows from Lemma \ref{quotient-fixedset} that $dA_i$ does not have $1$ 
as an eigenvalue. 
The map $\tatinvmp_i$, induced by $\tatinvmp$ on $G^{i-1}/G^i$, is given by 
$
\tatinvmp_i(tG^i)=\tatinv t G^i= (tG^i)A_i([tG^i]^{-1}) 
$. 
Since $G^{i-1}/G^i$ is abelian, $\tatinvmp_i$ is surjective. 

To prove the surjectivity of $\tatinvmp|_{G^{i-1}}$, let $b \in G^{i-1}$. 
Then there exists $w\in G^{i-1}$ with $bG^i=\tatinvmp_i(wG^i)=\tatinv w G^i$. This means $A(w)w^{-1}b$, and hence $w^{-1}bA(w)$ is an element of $G^i$. Now the surjectivity of $\tatinvmp|_{G^{i}}$ implies that there exists $y\in G^i$ such that $\tatinvmp(y)= \tatinv y=w^{-1}bA(w)$. Then we have 
\[
\tatinvmp(wy)=wyA(y^{-1})A(w^{-1})=ww^{-1}bA(w)A(w^{-1})=b.
\]
\end{proof}

Thus every element of the coset
$\coset A$ has finite order if $1$ is not an eigenvalue of $dA$. The asymptotic density of a coset is computed in the following corollary. 

\begin{cor}\label{equidistributed}
If $\coset A$ is any coset of $\Lambda$ in $\G$, then
\[
\lim_{r \ra \infty} \frac {\card{\coset A \cap \gammaball{S}{r}}}
{\card{\gammaball{S} {r}}}=
\frac{1}{\card F}.
\]
\end{cor}
\begin{proof}
Pick a set of coset representatives $\{\gamma_B\;|\;B\in F; \gamma_B\in \coset B\}$ and let $L=\max\{\ell_S(\gamma_B)\;|\;B\in F\}$. For any $A\in F$, there is a 
bijective map $\coset A \ra \Lambda$ given by $x \mapsto x\gamma_A^{-1}$. Then for any $r>0$, we have
\begin{equation}\label{coset-approx}
\card{\Lambda \cap \gammaball{S}{r-L}} \leq \card{\coset A \cap 
\gammaball{S}{r}} \leq \card{\Lambda \cap \gammaball{S}{r+L}}.
\end{equation}
Since $\card{\gammaball{S}{r}}= \sum_{A\in F}{\card{\coset A \cap \gammaball{S}{r}}}$, it follows that 
\begin{equation}\label{lambda-approx}
\card{\gammaball{S}{r-L}} \leq
\card F \card{\Lambda \cap \gammaball{S}{r}}\leq \card{\gammaball{S}{r+L}}. 
\end{equation}
A simple consequence of Theorem \ref{pansu} is that 
 $\displaystyle \lim_{r\ra \infty} 
\card{\gammaball{S}{r+ N}}/\card{\gammaball{S}{r}}=1$, for any 
$N \in \Z$. Now equation \eqref{lambda-approx} 
implies that $\displaystyle \lim_{r\ra \infty} 
\card{\Lambda \cap \gammaball{S}{r}}/\card{\gammaball{S}{r}}=1/{\card F}$ 
and the result follows from equation \eqref{coset-approx}. 
\end{proof}

Putting together the different pieces yields the formula for $\ratio \G S$:

\begin{proof}[End of proof of Theorem \ref{nilp}]
The inequality \eqref{cosetbound} was proved in Lemma \ref{estimate}. Now 
let $m$ be the number of elements of $\rho(F)$ which do not have $1$ as an eigenvalue.
Combining equation \eqref{contribution}, Lemma \ref{fullcoset}, and 
Corollary \ref{equidistributed}
we have 
\[
\ratio \G S
=\lim_{r \rightarrow \infty} \sum_{\substack {1 \mathrm{\, not \,an\, } \\ \text{eigenvalue}\\ \text{of } A}}
\frac{\card{\coset{A}\cap\gammafinball{S}{r}}}{\card{\gammaball{S}{r}}}
=m\lim_{r \ra \infty} 
\frac {\card{\coset A \cap \gammaball{S}{r}}}{\card{\gammaball{S} {r}}}
=\frac{m}{\card F}.
\]
\end{proof}

\section{Arbitrary virtually nilpotent groups}
\label{arbit}
Let $\G$ be any 
virtually nilpotent group. As discussed in Section \ref{dekigod}, $\Gamma$ has
a unique maximal finite normal subgroup, say $Q$,
and $\G/Q$ is almost crystallographic. So $\genfreeratio{\G/Q}$ can be computed using Theorem \ref{nilp}. Further, we have the following result.

\begin{cor}\label{allvnilpgps}
Let $\G$ and $Q$ be as above. Then 
$\ratio{\G}{S}=\genfreeratio {\G/Q}$ for any generating set $S$ of $\G$.
\end{cor}

\begin{proof}
Let $S= \{\g_1, \dots \g_l\}$ be a generating set for $\G$. Then 
$\bar S= \{\g_1Q, \dots \g_lQ\}$ generates $\G/Q$. 
Let $\ell_S$ and $\ell_{\bar S}$ denote the corresponding length functions on  $\G$ and $\G/Q$ respectively. 

Let $g\in \G$. Clearly,
$\ell_{\bar S}(gQ) \leq \ell_S(g)$. 
Moreover, 
if $\g_{i_1}Q\cdots \g_{i_n}Q$ is a geodesic word representing $gQ$, then $g=\g_{i_1}\cdots\g_{i_n}q'$ for some $q' \in Q$. If $M= \max \{ \ell_S(q)\;|\; q \in Q\}$, then $\ell_S(g) \leq  \ell_{\bar S}(gQ)+M$.

Let $\gammaball{\G}{r}$ and $\gammaball{\G/Q}{r}$ denote the balls of radius $r$ in $\G$ and $\G/Q$ respectively. 
Let $\gammafinball{\G}{r}$ and $\gammafinball{\G/Q}{r}$ represent the corresponding sets of finite-order elements. The above inequalities yield:
\[
\frac{|\gammaball{\G}{r}|}{|Q|}\leq |\gammaball{\G/Q}{r}| \text{ and }
|\gammaball{\G/Q}{r}|\leq \frac{|\gammaball{\G}{r+M}|}{|Q|}. 
\]
Since $Q$ is finite, an element of $\G$ has finite order if and only if its projection in $\G/Q$ has finite order. Thus we have
\[\frac{|\gammafinball{\G}{r}|}{|Q|}\leq |\gammafinball{\G/Q}{r}|
\text{ and }
|\gammafinball{\G/Q}{r}|\leq \frac{|\gammafinball{\G}{r+M}|}{|Q|}. 
\]
Putting together the above information, we have
\[
\frac{\card{\gammafinball{\G/Q}{r-M}}}{\card{\gammaball{\G/Q}{r}}}
\leq
\frac{\card{\gammafinball{\G}{r}}}{\card{\gammaball{\G}{r}}}
\leq
\frac{\card{\gammafinball{\G/Q}{r}}}{\card{\gammaball{\G/Q}{r-M}}}.
\]
Theorems \ref{nilp} and \ref{pansu} can now be used to conclude that
\[\ratio{\G}{S}=
\lim_{r\rightarrow  \infty}\frac{\card{\gammafinball{\G}{r}}}{\card{\gammaball{\G}{r}}}=\genfreeratio{\G/Q}.
\]
\end{proof}


\section{Examples}\label{examples}

\section{Examples}
}

Crystallographic groups form a large class of examples of virtually nilpotent groups.  
The results of the computation of $\genfreeratio \G$ for the crystallographic groups in dimensions $2$, $3$, and $4$, computed using GAP, are summarised in the Appendix.
The study of $\genfreeratio \G$ for crystallographic groups leads to a number of questions:
\begin{itemize}
\item Which numbers occur as $\genfreeratio \G$ for some crystallographic group $\G$?
\item What is the highest density that can occur in any dimension? 
\item What is the smallest dimension that a given density occurs in?
\item Is there an interesting explanation for the spectrum of densities in a given dimension?
\end{itemize}

We address the first of these and give a partial answer for the second. 

Theorem \ref{nilp} shows that $\genfreeratio{\G}$ is always a rational number in $[0,1)$. We now construct examples  
to show that in fact, every such number occurs as $\genfreeratio{\G}$ for some crystallographic group $\G$.

\subsection{Constructing examples}
The finite quotient $F=\G/\Lambda$ is called the \emph{holonomy group} of $\G$. 
The holonomy group of a crystallographic group can be realised as a finite subgroup of $Gl(n,\Z)$. On the other hand, if $F$ is a finite subgroup of $Gl(n,\Z)$, an averaging argument can be used to show that $F$ preserves an inner product on $\R^n$. Equivalently, there exists $M\in Gl(n,\R)$ such that 
$F'=MFM^{-1} \subset O(n)$. Then the lattice $\Lambda=M\Z^n$ is preserved by $F'$ and $\G=\Lambda \rtimes F'$ defines a crystallographic group.

The following Lemma is useful for constructing many examples. 

\begin{lem}\label{product}
Let $\G_1$ and $\G_2$ be virtually nilpotent groups. Then 
\[\genfreeratio{\G_1 \times \G_2}=\genfreeratio{\G_1}\genfreeratio{\G_2}.\]
\end{lem}

\begin{proof}In light of Corollary \ref{allvnilpgps}, we may assume $\G_i$ acts geometrically on a nilpotent Lie group $G_i$, for $i=1,2$. 
In this case $\G_1 \times \G_2$ acts geometrically on $G_1 \times G_2$.  
If $\G_i$ fits into 
\[
0 \rightarrow  \Lambda_i \rightarrow \G _i \rightarrow F_i \rightarrow 1,
\]
where $\Lambda_i$ is maximal normal nilpotent, then we have 
\[
0 \rightarrow  \Lambda_1 \times \Lambda_2
\rightarrow \G _1 \times \G _2 \rightarrow F_1 \times F_2 \rightarrow 1,
\]
and $\Lambda_1 \times \Lambda_2$ is the maximal normal nilpotent subgroup of $\G_1 \times \G_2$. 

If $A=(A_1,A_2) \in F_1 \times F_2$, then the set of eigenvalues of $dA$ is the union of the eigenvalues of $dA_1$ and $dA_2$. In particular, $1$ is not an eigenvalue for $dA$ if and only if neither $dA_1$ nor $dA_2$ has $1$ as an eigenvalue. Thus
$\genfreeratio {\G_1 \times \G_2}=\genfreeratio{\G_1}\genfreeratio{\G_2}$.
\end{proof}

\begin{proof}[Proof of Theorem \ref{allrationals}]
We start by constructing, for any $m\in \Z$, a crystallographic group $\G_m$
such that $\genfreeratio{\G_m}=\frac{m-1}{m}$. 
Let $\zeta$ be a primitive $m$th root of unity. 
If $\Phi$ denotes the Euler function, then $\{1, \zeta, \zeta^2 \dots \zeta^{\Phi(m)-1}\}$ is a basis for $\mathbb{Z}[\zeta]$, and we have 
$ \zeta^{\Phi(m)}= \sum_{i=0}^{\Phi(m)-1}a_i \zeta^i$, where $a_i \in \Z$. The matrix $T$, representing multiplication by $\zeta$ on $\Z[\zeta]$ is given below.

\[
T=
\begin{pmatrix}
 0&  &  &       &  &   a_0\\
 1&0 &  &       &  &   a_2\\
  &1 & 0 &       &  &  a_3\\
  &  & 1 &\ddots &  &   \vdots\\
  &  &  & \ddots &0 &   a_{\scriptscriptstyle \Phi(m)-2}\\
  &  &  &       &1 &a_{\scriptscriptstyle \Phi(m)-1}  \\
\end{pmatrix}
\]

The characteristic polynomial of $T$ is $ x^{\Phi(m)}-\sum_{i=0}^{\Phi(m)-1}a_ix^i$, which is also the minimal polynomial of $\zeta$. The eigenvalues of $T$ are exactly the $\Phi(m)$ primitive $m$th roots of unity. Thus $T$ has order $m$ and the matrices $T^i$ do not have $1$ as an eigenvalue, for $i<m$. 

Let $F\in O(n)$ be conjugate to $\langle T \rangle$ and let $\Lambda \simeq \Z^n$ be the lattice preserved by $F$. Then $\G_m=\Lambda \rtimes F$ is the desired group. 

Now let $\frac p q\in [0,1)$ and
note that $\frac p q = \frac{p}{p+1} \frac{p+1}{p+2} \cdots \frac{q-1}{q}$. Appealing to Lemma \ref{product} we can construct an example of a group $\G$ with 
$\genfreeratio{\G}=\frac p q$. 
\end{proof}

The above construction gives a very high dimensional crystallographic group if $p$ is much smaller than $q$. It would be nice to obtain a more efficient example. 

\subsection{Highest densities}
As seen in the tables in the appendix, the highest values of $\genfreeratio{\G}$ in two-, three- and four-dimensional crystallographic groups are $5/6$, $1/2$, and $23/24$ respectively.  The fact that the highest density in three dimensions is $1/2$ is part of a more general phenomenon:

\begin{prop}
If $\G$ is
 an odd-dimensional crystallographic group, $\genfreeratio \G \leq 1/2$.
\end{prop}
\begin{proof}
The holonomy group $F$ of $\G$ can be realised as a finite subgroup of $Gl(n,\Z)$. Since elements of $F$ have finite order, all their eigenvalues are roots of unity. Thus if
$n$ is odd and $A\in F\subset Gl(n,\Z)$ is orientation preserving (i.e. $A$ has determinant $1$), then $1$ is necessarily an eigenvalue of $A$. Thus at least half the elements of $F$ have $1$ as an eigenvalue, proving the result. 
\end{proof}

The upper bound is attained, for example by $\Z^n \rtimes \Z_2$, where the non-trivial element of $\Z_2$ is the automorphism $T$ of $\Z^n$ defined by $T(v)=-v$ for all $v$.  

\subsection{Non-abelian nilpotent groups}
We now investigate $\genfreeratio{\G}$ for some almost crystallographic groups 
$\G$. For these, the holonomy group can be realised as a finite group of automorphisms of the associated Lie algebra. We first show that in three and four dimensions, this turns out to be too restrictive, and $\genfreeratio{\G}$ is always $0$.


Recall that the \emph{complexification} of $\mathfrak g$, denoted by $\mathfrak g_{\mathbb C}$, is
$\mathfrak g \otimes_\R \mathbb C$. Any inner product on $\mathfrak g$ 
extends to a positive definite hermitian form on $\mathfrak g_{\mathbb C}$.  Any inner-product-preserving automorphism $T$ of $\mathfrak g$ extends to a unitary 
operator with the same eigenvalues. Further, if $\lambda$ is an eigenvalue of 
$T$, then there is an eigenvector (or a generalised eigenspace of the appropriate dimension) corresponding to $\lambda$ in $\mathfrak g_{\mathbb C}$.

\begin{lem}\label{eigensnorm1}
If $\mathfrak g$ is a nilpotent Lie algebra with inner product $\langle . \rangle$ and $T$ is an inner-product-preserving automorphism, then every eigenvalue of $T$ has norm $1$.
\end{lem}


\begin{proof}
Let $\lambda$ be an eigenvalue for $T$. Passing to the complexification of $\mathfrak g$, let $v$ be an eigenvector for $\lambda$. Then
$
\langle v,v\rangle =\langle Tv, Tv \rangle=\langle \lambda v, \lambda v\rangle= \lambda \bar \lambda\langle v,v\rangle
$
so that $|\lambda|=1$.
\end{proof}

\begin{lem}\label{1iseigen}
Let $T$ be an automorphism of a $3$- or $4$-dimensional nilpotent, non-abelian Lie algebra $\mathfrak g$. If all eigenvalues of $T$ have norm $1$, then $1$ is an eigenvalue of $T$.
\end{lem}

\begin{proof} In each case below, we pass to the complexification of $\mathfrak g $ to ensure the existence of eigenvectors or generalised eigenspaces. 

If $\mathfrak g$ is $3$-dimensional, it is isomorphic to the Heisenberg Lie algebra. If $Z$ generates $\mathfrak g^2$, then $T(Z)= \pm Z$. 
Thus we may assume that $T(Z)=-Z$ and that the eigenvalues of $T$ are $-1$, $\lambda_1$ and $\lambda_2$, where $\lambda_1$ and $\lambda_2$ are either both real, or complex conjugates of each other. 

If $\lambda_1$ and $\lambda_2$ are distinct, there exist distinct eigenvectors
$v_1$ and $v_2$ in $\mathfrak g_{\mathbb C}$. Then $[v_1,v_2]=cZ$ for some $c \in \mathbb C$. Since $T$ preserves the bracket, $[Tv_1,Tv_2]=T(cZ)=-cZ$. On the other hand, $[Tv_1,Tv_2]=[\lambda_1v_1,\lambda_2 v_2]=\lambda_1\lambda_2[v_1,v_2]=
\lambda_1\lambda_2 cZ$. We conclude that $\lambda_1\lambda_2=-1$. This cannot happen if $\lambda_1$ and $\lambda_2$ are conjugates, so the two eigenvalues 
have to be 
$1$ and $-1$.

If $\lambda_1=\lambda_2=\lambda$, then there exist vectors $v_1$ and $v_2$ in $\mathfrak g_{\mathbb C}$ such that $T(v_1)=\lambda v_1$ and 
$Tv_2=\lambda v_2+\alpha v_1$, for some $\alpha$. 
Let $[v_1,v_2]=cZ$ for some $c \in \mathbb C$. 
Then
$-cZ=[Tv_1,Tv_2]=[\lambda v_1,\lambda v_2+\alpha v_1]
=\lambda^2[v_1,v_2]=\lambda^2cZ$, which is impossible, since $\lambda$ is real. 

If $\mathfrak g$ is $4$-dimensional and two-step nilpotent, then $\mathfrak g^2$ is necessarily $1$-dimensional. If $Z$ generates $\mathfrak g^1$, then $T(Z) = \pm Z$. Assuming $T(Z)=-Z$, so that $-1$ is an eigenvalue, at least one of the other eigenvalues must be real. 
Thus we may assume the set of eigenvalues is $\{ -1, -1, \lambda_1, \lambda_2\}$ and proceed as above. 


If $\mathfrak g$ is $3$-step nilpotent, $\mathfrak g^3$ and $\mathfrak g^2/\mathfrak g^3$ are $1$-dimensional. Let $\mathfrak g^3=\langle Z \rangle$ and 
$\mathfrak g^1=\langle W ,Z \rangle$.  Then $T(Z)=\pm Z$ and $T(W)=\pm W +aZ$, for some $a \in \R$.  We may assume the set of eigenvalues is $\{-1, -1, \lambda_1, \lambda_2\}$. An argument similar to the above completes the proof. 
%
\end{proof}

\comment{
\begin{proof}
Once again, we pass to the complexification of the Lie algebra, to ensure the existence of eigenvectors. The bracket remains the same, and $A$ extends to an automorphism of the new Lie algebra, with the same eigenvalues. (think about this some more).
Pick a basis $\{X,Y,Z\}$ for $\mathfrak h$ with $[X,Y]=Z$. Since the central direction must be preserved, we can assume that $A(Z)=-Z$, i.e. that $-1$ is one of the eigenvalues. Suppose the other eigenvalues are $\lambda_1$ and $\lambda_2$. If they are distinct, we can find eigenvectors $v_1$ and $v_2$, corresponding to $\lambda_1$ and $\lambda_2$ respectively. 

Now $[v_1,v_2]=cZ$ for some $c$. So we should have $-cZ=[Av_1,Av_2]
=[\lambda_1 v_1,\lambda_2v_2]=\lambda_1\lambda_2[v_1,v_2]=\lambda_1\lambda_2cZ$. 
i.e. $\lambda_1\lambda_2=-1$. This is only possible if the eigenvalues are $1$ and $-1$. 

If the eigenvalues are the same, say $\lambda$, let $v_1$ be an eigenvector and $\{v_1,v_2\}$ a basis for the generalised eigenspace. Let $Av_2=\alpha v_1+\lambda v_2$.(make sure it always has to be of this form).
Then we have the following contradiction: 
$-cZ=[Av_1,Av_2]=[\lambda v_1,\alpha v_1 +\lambda v_2]
=\lambda^2[v_1,v_2]=\lambda^2cZ$. 

\end{proof}
}

The following corollary follows from Theorem \ref{nilp} and 
Lemmas \ref{eigensnorm1} and \ref{1iseigen}.

\begin{cor}
If $\G$ is a group acting geometrically on a $3$- or $4$-dimensional 
nilpotent (non-abelian) Lie group, then $\genfreeratio{\G}=0$.
\qed
\end{cor}

We now construct a class of almost crystallographic groups $\G$, such that 
$\genfreeratio \G$ is non-zero.

\begin{defn}{\bf(Generalised Heisenberg Lie algebras)}
Let $\mathfrak h_n$ be the Lie algebra generated by $\{X_1, \dots, X_n, Y_1,\dots,Y_n,Z\}$ such that  
$[X_i,Y_i]=Z$ for $1\leq i\leq n$, and all other brackets are $0$. 
\end{defn}

The following construction can be done for any generalised Heisenberg Lie algebra $\mathfrak h_{2n}$, of dimension $4n+1$. We give the construction for $\mathfrak h_2$:

Define an automorphism $T$ on $\mathfrak h_2$ by
\[
\begin{matrix}
X_1 \mapsto \;\;\,X_2 & Y_1 \mapsto -Y_2 & Z \mapsto -Z\\ 
X_2 \mapsto -X_1 & Y_2 \mapsto \;\;\,Y_1 &\\
\end{matrix}
\]
It is easy to check that $T$ preserves the bracket. The matrix of $T$ with respect to the basis
$\{X_1,X_2,Y_1,Y_2,Z\}$ is given by:
\[
\left(
\begin{matrix}
0&-1&  & & \\
1& \;\;\,0&  & & \\
 &  & \;\;\,0&1& \\
 &  &-1&0& \\
 &  &  & &-1
\end{matrix}
\right)
\]
Thus $T$ is an automorphism of order $4$ whose eigenvalues are $\pm i$, and  
$-1$.

Let $H_2$ be the connected, simply connected nilpotent Lie group corresponding to $\mathfrak h_2$. Let $\widetilde T$ be the automorphism of $H_2$ with $d \widetilde {T}=T$ and $N=\exp(\Z^5)$ be the lattice in $H_2$ preserved by $\widetilde{T}$. Then  
 $\G=N\rtimes \langle \tilde T \rangle$ is an almost-crystallographic group with $\genfreeratio{\G}=\frac 1 2$.

Note that for any group acting on a Lie group 
with $1$-dimensional centre, the maximum value of $\mathcal F$ is $\frac 1 2$, since the square of any automorphism fixes the central direction. 
The groups $\G$ defined above attain this maximum value.

\section*{Appendix}
The computations summarised below were done using the
computer algebra software GAP \cite{gap} and  the software package ``Cryst'',
which contains
libraries of $2$-, $3$- and $4$-dimensional crystallographic groups.

The first table gives the values of $\genfreeratio \G$ for all the 
$2$-dimensional crystallographic groups. See \cite{martin} for a description of the notation. 
In the tables summarising the results in dimensions $3$ and $4$,
$\mathcal N(q)$ denotes the number of groups $\G$ for which $\genfreeratio 
\G=q$.

\[
\begin{array}[t]{cc}
\multicolumn{2}{c}{\text{Dimension two}} \\
\begin{array}[t]{|c|c|}
\multicolumn{1}{c}{\G} & \multicolumn{1}{c}{\genfreeratio \G} \\
\hline \Wa{1} & 0 \\
\hline \W{1}{1} & 0 \\
\hline \W{1}{2} & 0\\
\hline \W{1}{3} & 0\\
\hline
\multicolumn{2}{c}{ } \\
\hline \Wa{6} & 5/6 \\
\hline \W{6}{1} &  5/12\\
\hline
\multicolumn{1}{c}{} & \multicolumn{1}{c}{} \\
\hline \Wa{3} &  2/3\\
\hline \W{3}{1} & 2/3 \\
\hline \W{3}{2} & 1/3 \\
\hline
\end{array}
&
\begin{array}[t]{|c|c|}
\multicolumn{1}{c}{\G} & \multicolumn{1}{c}{\genfreeratio \G} \\
\hline \Wa{4} & 3/4 \\
\hline \W{4}{1} & 3/8 \\
\hline \W{4}{2} & 3/8 \\
\hline
\multicolumn{2}{c}{ } \\
\hline \Wa{2} & 1/2 \\
\hline \W{2}{1} & 1/4 \\
\hline \W{2}{2} & 1/4 \\
\hline \W{2}{3} & 1/4 \\
\hline \W{2}{4} & 1/4 \\
\hline\end{array}
\end{array}
\hspace{1in}
\begin{array}[t]{c}
\multicolumn{1}{c}{\text{Dimension three}} \\
\begin{array}[t]{|c|c|}
\multicolumn{1}{c}{q} & 
\multicolumn{1}{c}{\num q}\\
\hline 0 & 113\\
\hline 1/8 & 28 \\
\hline 1/6 & 4 \\ 
\hline 3/16 & 20 \\ 
\hline 5/24 & 4 \\
\hline 1/4 & 30 \\ 
\hline 5/16 & 10 \\ 
\hline 1/3 & 1 \\ 
\hline 3/8 & 13 \\
\hline 5/12 & 2 \\ 
\hline 1/2 & 5 \\
\hline
\end{array}
\end{array}
\]

}
\begin{table*}
\begin{center}
{Dimension four}
\label{fourdim}
$
\begin{array}[t]{|c|c|}
\multicolumn{1}{c}{q} & 
\multicolumn{1}{c}{\num q}\\
 \hline 0& 1875 \\ \hline 1/16& 605 \\ \hline 1/12& 64 \\ \hline 3/32& 426 \\ \hline 5/48& 48 \\
  \hline 1/9& 25 \\ \hline 1/8& 558 \\ \hline 5/36& 5 \\ \hline 9/64& 50 \\ \hline 5/32& 193 \\ 
  \hline 1/6& 38 \\ \hline 25/144& 2 \\ \hline 13/72& 7 \\ \hline 3/16& 229 \\ \hline 1/5& 10 \\ 
  \hline 13/64& 23 \\ \hline 5/24& 31 \\ \hline 2/9& 31 \\ \hline 9/40& 2 \\ \hline 15/64& 11 \\ 
  \hline 1/4& 125 \\ \hline
   \end{array}   
  \quad
  %
\begin{array}[t]{|c|c|}
\multicolumn{1}{c}{q} & 
\multicolumn{1}{c}{\num q}\\
 \hline 37/144& 3 \\ \hline 33/128& 16 \\ \hline 35/128& 8 \\ \hline 5/18& 9 \\ 
  \hline 9/32& 34 \\ \hline 23/80& 2 \\ \hline 85/288& 1 \\ \hline 5/16& 91 \\ \hline 21/64& 20 \\ 
  \hline 1/3& 20 \\ \hline 385/1152& 1 \\ \hline 49/144& 2 \\ \hline 25/72& 2 \\ \hline 205/576& 2 \\ 
  \hline 13/36& 6 \\ \hline 3/8& 27 \\ \hline 2/5& 6 \\ \hline 13/32& 9 \\ \hline 5/12& 4 \\ 
  \hline 7/16& 2 \\ \hline 4/9& 12 \\ \hline
   \end{array}   
\quad
%
\begin{array}[t]{|c|c|}
\multicolumn{1}{c}{q} & 
\multicolumn{1}{c}{\num q}\\
 \hline 9/20& 1 \\ \hline 11/24& 1 \\ \hline 15/32& 4 \\ \hline 1/2& 6 \\ 
  \hline 37/72& 2 \\ \hline 33/64& 10 \\ \hline 25/48& 1 \\ \hline 17/32& 7 \\ \hline 35/64& 5 \\ 
  \hline 5/9& 7 \\ \hline 9/16& 17 \\ \hline 23/40& 2 \\ \hline 85/144& 1 \\ \hline 43/72& 1 \\ 
  \hline 5/8& 14 \\ \hline 91/144& 1 \\ \hline 21/32& 13 \\ \hline 2/3& 4 \\ \hline 385/576& 1 \\ 
  \hline 65/96& 2 \\ \hline 49/72& 1 \\ \hline
   \end{array}   
\quad
%
\begin{array}[t]{|c|c|}
\multicolumn{1}{c}{q} & 
\multicolumn{1}{c}{\num q}\\
 \hline 11/16& 4 \\ \hline 25/36& 3 \\ \hline 17/24& 1 \\ \hline 205/288& 1 \\ 
  \hline 137/192& 2 \\ \hline 13/18& 2 \\ \hline 35/48& 3 \\ \hline 3/4& 3 \\ \hline 55/72& 2 \\ 
  \hline 19/24& 1 \\ \hline 4/5& 1 \\ \hline 77/96& 4 \\ \hline 13/16& 4 \\ \hline 5/6& 1 \\ 
  \hline 41/48& 3 \\ \hline 31/36& 1 \\ \hline 7/8& 4 \\ \hline 9/10& 1 \\ \hline 11/12& 2 \\ 
  \hline 23/24& 4 \\ \hline
   \end{array} 
$     
\end{center}
\end{table*}

%
%



\end{document}